\documentclass[12pt]{article} \usepackage{amssymb}
 \def\vbar{\mathchoice{\vrule height2.3ptdepth-.3ptwidth.12pt\kern-
 .10pt}
    {\vrule height6.3ptdepth-.3ptwidth.11pt\kern-.11pt}
    {\vrule height5.1ptdepth-.30ptwidth.8pt\kern-.8pt}
    {\vrule height4.1ptdepth-.24ptwidth.6pt\kern-.7pt}}
\newenvironment{pff}{\hspace*{-\parindent}{\bf Proof:}}{\hfill $\Box$
\vspace*{0.2cm}}

\def\reel{\hbox{I\hskip-2pt R}}
\def\comp{\hbox{I\hskip-6pt C}}

\def\tr{\;\mathrm{tr}\;}
\def\<{\langle}
\def\>{\rangle}

\def\n{{\boldmath n}}

\def\<{\langle}
\def\>{\rangle}

\def\reel{\hbox{I\hskip -2pt R}}

\def\tr{\textmd{trace}\,}

\def\n{{\noindent}}

\newtheorem{theorem}{Theorem}[section]

\newtheorem{lemma}{Lemma}[section]
\newtheorem{prop}{Proposition}[section]

\newtheorem{example}{Example}[section]

\newenvironment{Pff}{\hspace*{-\parindent}{\bf Proof}}{\hfill $\Box$
\vspace*{0.2cm}}

\def\tr{\;\mathrm{tr}\;}
\def\opp{\;\mathrm{opp}\;}

\begin{document}

 \vspace{0.5cm}
\title{{\bf
Extension of the Olkin and Rubin Characterization to the Wishart
distribution on homogeneous cones}}
\author{I. Boutouria$^*$, A. Hassairi\thanks{Sfax University,
Tunisia} \footnote{Corresponding author.
 \textit{E-mail address: Abdelhamid.Hassairi@fss.rnu.tn}}
 and H. Massam \thanks{Department of Mathematics and Statistics, York
University, Canada. This author was supported by NSERC grant A
8947.} }

\date{}
\maketitle{}

\begin{abstract}

The Wishart distribution on an homogeneous cone is a generalization
of the Riesz distribution on a symmetric cone which corresponds to a
given graph. The paper extends to this distribution, the famous
Olkin and Rubin characterization of the ordinary Wishart
distribution on symmetric matrices.
\end{abstract}

\n {\small {\it{ Keywords:}} Vinberg algebra, homogeneous cone,
Wishart distribution, orthogonal group.

\section{Introduction}

In many practical situations, there are manifest inter-relationships
among several variables. One important case is when several pair
variables are conditionally independent, giving other remaining
variables.  For multivariate normal distribution, this corresponds
to some zeros among the entries of the precision matrix. Due to this
there has been an interest in distributions akin to the Wishart but
defined more generally on various cones containing the cone $\Omega$
of positive definite symmetric matrices as a special case. In
particular, Andersson and Wojnar $[2]$ defined the Wishart
distribution on homogeneous cones. This distribution is in fact an
extension to homogeneous cone of the Riesz distribution on symmetric
cones defined by Hassairi and Lajmi in $[8]$. In the present paper,
we give a characterization of the Wishart distribution on
homogeneous cones which is parallel to that given in $[13]$ by Olkin
and Rubin, or more generally in $[6]$ by Casalis and Letac. In these
papers, it is not assumed that densities exist, however
distributions are assumed to be invariant by the orthogonal group of
the appropriate algebra. Our characterization uses the Laplace
transform and a decomposition of a random variable on an homogeneous
cone as a sum of random variables concentrated on certain
subalgebras. The distribution of each component is then assumed to
be invariant by the orthogonal group of the corresponding
subalgebra. Let us give, a brief history of this characterization.
Luckacs (see $[11]$) gave the following characterization of the
Gamma distribution: Let $X$ and $Y$ be two non Dirac and non
negative independent random variables
 such that $X+Y$ is positive almost surely, then $U=X+Y$ is independent of $V=\frac{X}{X+Y}$ if and only if there exist $\sigma>0$, $p>0$ and $q>0$ such that
$X$ and $Y$ are distributed as Gamma distributions with parameters
$(p,\sigma)$ and $(q,\sigma)$ respectively. Olkin and Rubin $[13]$
extended this characterization to the Wishart distribution on the
cone $\Omega$ of positive definite symmetric matrices. They showed
that $X$ and $Y$ in $\Omega$ are Wishart if and only if
$(X+Y)^{-1/2}X (X+Y)^{-1/2}$  is independent of $X+Y$ and its
distribution is invariant by the orthogonal group $K$. The
characterization has then been extended by Carter $[5]$ to the
Wishart on the cone of Hermitian matrices with entries in $\comp$,
and by Casalis and Letac $[6]$ to any symmetric cone. There are
other types of characterizations of the Wishart such as, for
example, that given by Letac and Massam in $[9]$, Geiger and
Heckerman in $[7]$, and Massam and Weso\l owski in $[12]$. We also
mention that more recently Bobecka and Weso\l owski gave in $[3]$ a
characterization of the Wishart distribution on $\Omega$ without any
assumption of invariance. However, they had to assume that the
densities of both $X$ and $Y$ with respect to the Lebesgue measure
exist and are twice differentiable. The extension of this
characterization to homogeneous cones has been the object of a paper
by Boutouria $[4]$. The paper is organized as follows. In \S 2, we
recall various definitions and preliminary results relevant to the
Wishart on homogeneous cones and we establish some results
concerning the Vinberg multiplication and determinant calculation.
In \S 3, we state and prove our main characterization result.

\section{The Wishart distribution on homogeneous cones}
For the convenience of the reader, we will give here the elements of
Vinberg algebras and homogeneous cones essential to working with the
family of Wishart distributions on homogeneous cones. These elements
are taken from $[2]$ and the reader is refereed to this paper for
further details. After recalling this development, we give three
examples. For the third example, the poset is isomorphic to a rooted
tree while this is not so for the second example. The first example
corresponds to the ordinary Wishart distribution on symmetric
matrices.

Let $I$ be a partially ordered finite set (herewith abbreviated as
poset) equipped with a relation denoted $\preceq$. We will write
$i\prec j$ if $i\preceq j$ and $i \not = j$. We assume that $I$
satisfies the following condition
\begin{eqnarray*}
(F): \left \{ \begin{array}{c}
 \mbox{for any two points}\;\; i\;\;\mbox{and}\;\; j\;\;\mbox{in}\;\; I \;\;\mbox{such that either}\;\;i\prec j\;\;\mbox{or}\;\;j\prec i\;\; \\
 \mbox{the path on the Hasse diagram
of}\;\;I\;\;\mbox{between}\;\; i\;\;\mbox{and}\;\; j\;\;\mbox{is
unique}.
\end{array}\right.
\end{eqnarray*}
For all pairs $(i,j)\in I \times I$ with $j\prec i$, let $E_{ij}$ be
a finite-dimensional vector space over $\reel$ with
$n_{ij}=\mbox{dim}(E_{ij})>0$. Set
\begin{eqnarray*}
{\mathcal{A}}_{ij}= \left \{ \begin{array}{ccccc}
  {\reel} & \textrm{for} & i=j&& \\
  E_{ij} & \textrm{for} & j\prec i\\
  E_{ji} & \textrm{for} & i\prec j \\
  \{0\} &  &\textrm{otherwise}&&
\end{array}\right.
\end{eqnarray*}
\n and ${\mathcal{A}}=\displaystyle\prod\limits _{i,j \in I\times
I}{\mathcal{A}}_{ij}$. \n Define $n_{i.}=\displaystyle\sum \limits
_{\mu \prec i}n_{i\mu}$, $n_{.i}=\displaystyle\sum \limits _{i\prec
\mu }n_{\mu i}, \ n_{i}=1+\frac{1}{2}(n_{i.}+ n_{.i})$, $i\in I $
and $n_{..}=\displaystyle\sum \limits _{i \in I}n_{i}$.

\n An element $A\equiv (a_{ij}, \ i,j \in I)$ of ${\mathcal{A}}$ may
be seen as a matrix and so we define its trace as $\textrm{tr}
A=\displaystyle \sum _{i \in I}a_{ii}$.

\n Let $f_{ij}\;:\; E_{ij} \rightarrow E_{ij}$, $i\succ j$, be
 involutional linear mappings, i.e., $f_{ij}^{-1}=f_{ij}$. They
induce an involutional mapping ( $ A \mapsto A^{\ast}$ ) of
${\mathcal{A}}$ given as follows: $A^{\ast}=(a^*_{ij}|(i,j)\in I
\times I)$, where

$$a^{\ast}_{ij}=\left \{\begin{array}{ccccc}a_{ii} \hskip1.5cm& \textrm{for} & i=j&&  \\
f_{ij}(a_{ji}) \hskip1cm& \textrm{for} &j\prec i\\
 f_{ji}(a_{ji})\hskip1cm& \textrm{for}&i\prec j \\
  0\hskip1.5cm&   & \hskip-0.8cm\textrm{otherwise.}&& \end{array}\right.$$

\n We now define the following subspaces of ${\mathcal{A}}$:

\n the upper triangular matrices
\begin{eqnarray}\label{r1}
 {\mathcal{T}}_{u}=\{A\equiv(a_{ij})\in{\mathcal{A}},\  \forall i,j \in I :i
\not\preceq j \Rightarrow a_{ij}=0\};
\end{eqnarray}
\n the lower triangular matrices
\begin{eqnarray}\label{r2}
{\mathcal{T}}_{l}=\{A\equiv(a_{ij})\in{\mathcal{A}},\   \forall i,j
\in I :j\not\preceq i  \Rightarrow a_{ij}=0\}; \
\end{eqnarray}
\n and the Hermitian matrices ${\mathcal{H}}=\{A \in{\mathcal{A}}, \
A^{\ast}=A\}. $

\n The sets of upper and lower triangular matrices in
${\mathcal{A}}$ with positive diagonal elements are respectively
denoted by ${\mathcal{T}}_{u}^{+}$ and ${\mathcal{T}}_{l}^{+}$. The
sets of upper and lower triangular matrices with all diagonal
elements equal to 1 are respectively denoted by
${\mathcal{T}}_{u}^{1}$ and ${\mathcal{T}}_{l}^{1}$. The sets of
diagonal matrices and of diagonal matrices with positive entries are
denoted by ${\mathcal{D}}$ and ${\mathcal{D}}^{+}$, respectively.

\n We are going to equip the vector space ${\mathcal{A}}$ with a
bilinear map called multiplication and denoted by $(A,B)\mapsto AB$.
For this purpose we need to define bilinear mappings
${\mathcal{A}}_{ij}\times {\mathcal{A}}_{jk}\rightarrow
{\mathcal{A}}_{ik}$, denoted by $(a_{ij},b_{jk})\mapsto
a_{ij}b_{jk}, $ and then define $AB=C\equiv (c_{ij}|(i,j)\in I\times
I)$ by $c_{ij}=\displaystyle\sum \limits _{\mu \in I}a_{i\mu}b_{\mu
j}$.

\n The multiplication is required to satisfy the following
properties:
\begin{eqnarray}
\label{vinberg} &&i)\ \  \forall A\in {\mathcal{A}}; \ A \not= 0
\Rightarrow
\textrm{tr}(AA^{\ast})>0\nonumber\\
&&ii) \ \  \forall A, B \in {\mathcal{A}}; \
(AB)^{\ast}=B^{\ast}A^{\ast}\nonumber\\
&&iii)\ \   \forall A, B \in {\mathcal{A}}; \ \textrm{tr} (AB)=\textrm{tr}(BA)\\
&&iv) \ \  \forall A, B, C \in {\mathcal{A}}; \  \textrm{tr} (A(BC))=\textrm{tr} ((AB)C)\nonumber\\
&&v)\ \   \forall U,S,T \in {\mathcal{T}}_{l};  \ (ST)U=S(TU)\nonumber\\
&&vi) \ \  \forall U,T \in {\mathcal{T}}_{l}; \
T(UU^{\ast})=(TU)U^{\ast}.\nonumber
\end{eqnarray}
\n An algebra ${\mathcal{A}}$ with the above structure and
properties is called a Vinberg algebra (For more details, we can
refer to $[2]$). We choose the element $A\equiv(a_{ij}|(i,j)\in
I\times I)$ of ${\mathcal{D}}$ such that $a_{ii}=1, \ \forall \ i
\in I$ as the unit element of ${\mathcal{A}}$ and we denote it by
$e$. Vinberg proved in $[14]$ that the subset
${\mathcal{P}}=\{TT^{\ast}\in {\mathcal{A}}, \ T \in
{\mathcal{T}}_{l}^{+}\}$ $\subset \mathcal{H} \subset \mathcal{A}$
forms a homogeneous cone, that is the action of its automorphism
group is transitive.

\n The definition of ${\mathcal{P}}$ could be changed to the
following equivalent definition
$${\mathcal{P}}=\{TDT^{\ast} \in {\mathcal{A}},  \ T \in T_{l}^{1}, \  D \in {\mathcal D}^{+}\}. $$
\n The two decompositions $S=TT^{\ast}, \ T
\in{\mathcal{T}}_{l}^{+}\ $ and $S=T_{1}DT_{1}^{\ast}, \ T_{1} \in
{\mathcal{T}}^{1}_{l}, \ D \in {\mathcal{D}}^{+}$ are unique and
their connection is given by $T=T_{1}\sqrt{D}$ where $\sqrt{D}\equiv
\mbox{diag}(\sqrt{d_{i}}, \ i \in I) \in {\mathcal{D}}^{+}$ when $D
\equiv \mbox{diag}(d_{i}, \ i \in I) \in {\mathcal{D}}^{+}$.

\n For $S=(s_{ij}, \ i, j \in I)=T_{1}DT_{1}^{\ast}$, we write
$D_{ii}=S_{[i].}$.

If $\preceq^{\opp}$ is the opposite ordering on the index set $I$,
i.e., $i\preceq^{\opp}j\Leftrightarrow j\preceq i$. The Vinberg
algebra ${\cal{A}}^{\opp}=\displaystyle\prod\limits _{i,j \in
I\times I}{\mathcal{A}}^{\opp}_{ij}$, where
\begin{eqnarray*}
{\mathcal{A}}^{\opp}_{ij}= \left \{ \begin{array}{ccccc}
  {\reel} & \textrm{for} & i=j&& \\
  E_{ji} & \textrm{for} & j\prec^{\opp} i & \\
  E_{ij} & \textrm{for} &i\prec^{\opp} j\\
  \{0\} &  &\textrm{otherwise,}&&
\end{array}\right.
\end{eqnarray*}
differs from the Vinberg algebra ${\mathcal{A}}$ only in the
ordering of the index set $I$. It is proved (see $[12]$) that
${\mathcal{P}}_{\preceq^{\opp}}=\{T^{\ast}T\in {\mathcal{A}}, \ T
\in {\mathcal{T}}_{l}^{+}\}$ is the dual cone of ${\cal{P}}$ which
is also denoted ${\cal{P}}^*$.

Let's explain why condition (F) on the poset $I$ is required for the
definition of a Vinberg algebra. In fact, the property $vi)$ in
(\ref{vinberg}) fails to be verified if the condition (F) is not
satisfied. Suppose that (F) is not satisfied, then there exist two
elements $i$ and $j$ in $I$ such that either $i\prec j$ or $j\prec
i$ and the path on the Hasse diagram of $I$ between $i$ and $j$ is
not unique. Without loss of generality, we suppose that $i\prec j$.
Then there exist $k$ and $s$ in $I$ such that $k \not=s$, $i\prec
k\prec j$ and $i\prec s \prec j$. Consider the elements
$T=(t_{nm})_{n,m \in I}$ and $U=(u_{nm})_{n,m \in I}$ of
${\cal{T}}_l$ defined by $t_{nm}\not=0$ and $u_{nm}\not=0$ if $n,m
\in \{i,j,k\}$ and $t_{nm}=0$ and $u_{nm}=0$ otherwise. Then
$$[T(UU^*)]_{jk}=t_{ji}u_{ii}u_{ik}+t_{jk}(u_{ik}^2+u_{kk}^2)+t_{jj}(u_{ki}u_{ij}+u_{jk}u_{kk})$$
and
$$[(TU)U^*]_{jk}=(t_{ji}u_{ii}+t_{jk}u_{ki}+t_{js}u_{si}+t_{jj}u_{ji})u_{ik}+(t_{jk}u_{kk}+t_{jj}u_{jk})u_{kk}.$$
We see that $T(UU^*)\not =(TU)U^*$ so that $vi)$ in (\ref{vinberg})
is not satisfied.

Andersson, Letac and Massam $[1]$ have also shown, without using the
theory of Vinberg algebras, that if (F) is satisfied then, the cone
${\cal P}$ is homogeneous.

Let $G$ be the connected component of the identity in
Aut$({\cal{P}})$, the group of linear transformations leaving ${\cal
P}$ invariant. We recall that $\chi: G\mapsto \reel_+$ is said to be
a multiplier on the group $G$ if it is continuous, $\chi(e)=1$ and
$\chi(h_1h_2)=\chi(h_1)\chi(h_2)$ for all $h_1,h_2 \in G$. Consider
the map $\pi\;:\; T\in {\mathcal T}_{l}^+ \mapsto \pi(T)\in
\pi({\mathcal T}_{l}^+)\subset G$ such that for $X=WW^{\ast}\in
{\mathcal{P}}$, $W\in {\mathcal T}_{l}^+$
\begin{eqnarray}\label{o6}
\pi(T)(X)=(TW)(W^*T^*).
\end{eqnarray}

\n  Andersson and Wojnar show in $[2]$ that the restriction of a
multiplier $\chi$ to the (lower) triangular group ${\mathcal
T}_{l}^+$, i.e., $\chi \circ \pi : T_{l}^{+} \rightarrow
{\reel}_{+}$ is in one to one correspondence with the set of
$(\lambda_{i}, i \in I)\in {\reel}^{I}$. To each $\chi$ corresponds
a unique (up to a multiplicative constant) equivariant measure
$\nu^{\chi}$ on $C$ with multiplier $\chi$ under the action of $G$,
that is for all $h \in G$, the image measure $h^{-1} \nu^{\chi}$ of
$\nu^{\chi}$ by $h^{-1}$ is
$$h^{-1} \nu^{\chi}=\chi(h)\nu^\chi.$$

\n For $\theta\in {\cal P}^*$, the Laplace transform of $\nu^{\chi}$
is
\begin{equation}\label{a2}
L^{\chi}(\theta)=\displaystyle\int_{\cal
P}\exp\{-\theta(P)\}d\nu^{\chi}(P).
\end{equation}

\n We define the $\chi$-inverse of $\theta$ by
$$\theta^{\chi}=-\frac{d}{d \theta}\log L^{\chi}(\theta).$$
The mapping
\begin{eqnarray*}
\theta\in {\cal P}^*&\mapsto&\theta^{\chi}\in {\cal P},
\end{eqnarray*}
is a bijection. Its inverse  is denoted by
\begin{eqnarray}
P\in {\cal P}&\mapsto& P^{-\chi}\in {\cal P}^*,\nonumber
\end{eqnarray}
\n and we have the properties
\begin{equation}\label{c1}
(h^{-1}P)^{-\chi} = \ ^{t}hP^{-\chi}\ \ \textrm{and}
 \ \ (^{t}h\theta)^{\chi}=h^{-1}(\theta^{\chi}),
\end{equation}
where $h \in G$, $P \in {\cal {P}}$ and $\theta \in {\cal {P}}^{\ast}.$
Also for $\theta\in {\cal P}^*$, it is convenient to introduce the
notation $\sigma=\theta^{-\chi}$. Then (\ref{a2}) may be written as
a function of $\sigma$ as
\begin{eqnarray}
n^{\chi}(\sigma)&=&\int_{\cal
P}\exp\{-\sigma^{-\chi}(P)\}d\nu^{\chi}(P),\nonumber
\end{eqnarray}
so that for $h \in G, \sigma \in {\cal P}$,  $n^{\chi}$ has the
property
\begin{eqnarray}\label{o7}
n^{\chi}(h\sigma)=\chi(h)n^{\chi}(\sigma).
\end{eqnarray}
Andersson and Wojnar (see $[2]$) consider the set
$${\mathcal{X}}=\{\chi \circ \pi :{\mathcal{T}}_{l}^{+} \rightarrow {\reel}_{+};  \
\lambda_{i}>\frac{n_{i.}}{2}, \ i \in I\},$$ and show that for $\chi
\in {\mathcal{X}}$, the measure concentrated on  ${\mathcal{P}}$
\begin{equation}
\label{meas} \nu^{\chi}(dX)=\prod\limits_{i\in
I}x_{[i].}^{\lambda_{i}-n_{i}}{ 1}_{\mathcal{P}}(X)dX
\end{equation}
generates the Wishart natural exponential family of distributions on
${\mathcal{P}}$ absolutely continuous with respect to the Lebesgue
measure, parameterized by $\sigma \in {\mathcal{P}}$
\begin{equation}
\label{wishart} HW_{\chi,\sigma}(dX)=\frac{\pi^{\frac{|I|-n_{\cdot
\cdot}}{2}}\prod\limits_{i\in
I}\lambda_{i}^{\lambda_{i}}\prod\limits_{i\in I}
x_{[i].}^{\lambda_{i}-n_{i}}}{\prod \limits_{i\in
I}\Gamma(\lambda_{i}-\frac{n_{i.}}{2})\prod\limits_{i\in I}
\sigma_{[i].}^{\lambda_{i}}} \exp {-\{{\tr (\sigma^{-\chi} }X )\}}{
\textbf{1}}_{\mathcal{P}}(X)dX.
\end{equation}
It is shown in $[2]$ that if $\sigma=ZZ^\ast$, where $Z\in {\cal
T}_l^{+}$,
\begin{equation}\label{sigma-chi}
\label{sigma} \sigma^{-\chi}=(Z^\ast)^{-1}\mbox{diag}(\lambda_i|i\in
I)Z^{-1}.
\end{equation}
We are now going to give three examples of homogeneous cones and
their corresponding Wishart distributions.
\begin{example}\label{special case}
 In this example, we show how the ordinary Wishart distribution on symmetric
matrices is a particular case of the general Wishart distribution
$HW_{\chi,\sigma}$ on homogeneous cones. Take
$I=\{1,\cdot\cdot\cdot,I\}$ equipped with the usual total ordering
$\preceq$ on integers, denote by $I$ its cardinality, and set
$E_{ij}=\reel, \ j\prec i$. Then the vector space ${\cal{A}}$ is the
space ${\cal{M}}(I,\reel)$ of all $I\times I$ matrices. With the
standard multiplication and inner product the vector space
${\cal{M}}(I,\reel)$ is a Vinberg algebra. The homogeneous cone
${\cal{P}}$ in this Vinberg algebra is the usual cone of $I \times
I$ positive definite symmetric matrices usually denoted $\Omega$. In
this case, we have $$n_{.i}=I-i, \ n_{i.}=i-1, \ n_i=\frac{I+1}{2},
\ n_{..}=\frac{I(I+1)}{2}$$ and every multiplier $\chi: \ G
\rightarrow \reel_{+}$ has the form $\chi(h)=|\det(h)|^{\alpha}$,
$\alpha \in \reel$. Thus for all $i\in I$, $\lambda_{i}$ doesn't
depend on $i$, it is equal to $\lambda=\frac{I+1}{2}\alpha$, so that
$\chi$ can be replaced by $\lambda$. From (\ref{sigma-chi}), it
follows that $\sigma^{-\chi}=\lambda\sigma^{-1}$, where
$\sigma^{-1}$ is the usual inverse of $\sigma$ in ${\cal{P}}$. Since
for $S=(s_{ij}, \ \ i,j \in I) \in {\Omega}$,
$\det(S)=\displaystyle\prod_{i \in I}s_{[i].}$, the density of the
Wishart distribution (\ref{wishart}) becomes
\begin{eqnarray}\label{wishart usuel}W_{\lambda, \
\sigma}(dS)=\frac{\lambda^{I\lambda}\det(S)^{\lambda-\frac{I+1}{2}}}{\pi^\frac{I(I-1)}{4}\displaystyle\prod_{i
\in
I}\Gamma(\lambda-\frac{i-1}{2})\det(\sigma)^{\lambda}}\exp\{-\lambda\mbox{tr}(\sigma^{-1})\}{
\textbf{1}}_{\Omega}(S)dS,
\end{eqnarray}which is the usual
Wishart distribution on $\Omega$ with parameters
$\lambda>\frac{I-1}{2}$ and $\sigma \in \Omega$.
\end{example}
\begin{example}
Let $I=\{1,2,3,4\}$ be an index set and consider the poset defined
by
$$1\prec 3, 2 \prec 3, 2\prec 4.$$
The homogeneous cone corresponding to this poset cannot be of the
type $Q_G$ of incomplete matrices with submatrices corresponding to
the cliques of $G$ being positive definite. We will not give the
details of the argument for this assertion here. However, an
accurate, albeit short, argument is that according to Theorem 2.2 in
$[10]$, this can only happen if the undirected graph obtained by
dropping the directions on the graphical representation of the poset
is a homogeneous graph. The graphical representation of the poset in
this example is the directed graph with directed edges
$$\{(1,3), (2,3), (2,4)\}.$$
The undirected graph obtained by dropping the directions is the
three-link chain with undirected edges $\{(1,3), (3,2), (2,4)\},$
which, by definition, is not a homogeneous graph.

Let us illustrate this in the case where $E_{ij}=\reel$ for all
appropriate $(i,j)$. The Vinberg algebra $\cal{A}$ is thus the
vector space of all $I\times I$ real matrices with zeros at the
(1,2), (2,1), (3,4) and (4,3) entries. An element $X=(x_{ij})_{i,j
\in I}$ of ${\cal P}$ is decomposed into $X=TDT^\ast$ with $T\in
{\cal T}_l^1, D\in {\cal D}$, such that
$$
T^\ast=\pmatrix{1&0&t_{13}&0\cr 0&1&t_{23}&t_{24}\cr 0&0&1&0\cr
0&0&0&1}, D=\pmatrix{d_1&0&0&0\cr 0&d_2&0&0\cr 0&0&d_3&0\cr
0&0&0&d_4}.$$
 Then the equation
\begin{eqnarray}
X&=&\pmatrix{x_1&0&x_{13}&0\cr 0&x_2&x_{23}&x_{24}\cr
x_{31}&x_{32}&x_3&0\cr 0&x_{42}&0&x_4}=\pmatrix{1&0&0&0\cr
0&1&0&0\cr t_{31}&t_{32}&1&0\cr 0&t_{42}&0&1} \pmatrix{d_1&0&0&0\cr
0&d_2&0&0\cr 0&0&d_3&0\cr 0&0&0&d_4}\pmatrix{1&0&t_{13}&0\cr
0&1&t_{23}&t_{24}\cr 0&0&1&0\cr 0&0&0&1}
\nonumber\\
&=&\pmatrix{d_1&0&d_1t_{13}&0\cr 0&d_2&d_2t_{23}&d_2t_{24}\cr
t_{31}d_1&t_{32}d_2&t_{32}d_2t_{23}+d_3&t_{32}d_2t_{24}\cr
0&t_{42}d_2&t_{32}d_2t_{23}+d_3&t_{32}d_2t_{24}\cr
0&t_{42}d_2&t_{42}d_2t_{23}&t_{42}d_2t_{24}+d_4}\nonumber
\end{eqnarray}
has a unique solution
\begin{eqnarray}
d_1=x_1=x_{[1]\cdot}, \  \ t_{13}=\frac{x_{13}}{x_1},&&\;d_2=x_2=x_{[2]\cdot},\; t_{23}=\frac{x_{23}}{x_2},\; t_{24}=\frac{x_{24}}{x_2},\; \nonumber\\
\label{cdot} d_3=x_3-x_{32}x_2^{-1}x_{23}=x_{[3]\cdot},&&
d_4=x_4-x_{42}x_2^{-1}x_{24}=x_{[4]\cdot}.
\end{eqnarray}
The cone ${\cal{P}}$ is then the set of symmetric matrices
$X=(x_{ij})$ such that
\begin{eqnarray*}
x_1>0, \; x_2>0, \; x_3-x_2^{-1}x_{23}^2>0 \ \textrm{and} \
x_4-x_2^{-1}x_{24}^2>0.
\end{eqnarray*}
Since $E_{ij}=\reel$, we have that $n_{ij}=1$, $n_{1\cdot}=0, \
n_{2\cdot}=0, \ n_{3\cdot}=2$ and $n_{4\cdot}= 1.$ Therefore
$${\mathcal{X}}=\{(\lambda_1,\lambda_2,\lambda_3,\lambda_4):\; \lambda_1>0,\lambda_2>0,\lambda_3>1, \lambda_4>\frac{1}{2}\}.$$
 and the Wishart distribution on this homogeneous cone
 is given for $\chi \in {\mathcal{X}}$ and $\sigma \in {\cal{P}}$ by
\begin{eqnarray*}
HW_{\chi,\sigma}(dX)=\frac{\pi^{\frac{4-8}{2}}\prod\limits_{i=1
}^{4}\lambda_{i}^{\lambda_{i}}
x_{[1].}^{\lambda_{1}-\frac{3}{2}}x_{[2].}^{\lambda_{2}-2}x_{[3].}^{\lambda_{3}-2}
        x_{[4].}^{\lambda_{4}-\frac{3}{2}}  }{\Gamma(\lambda_{1})\Gamma(\lambda_{2})\Gamma(\lambda_{3}-1)\Gamma(\lambda_{4}-\frac{1}{2})
 \prod\limits_{i=1}^{4} \sigma_{[i].}^{\lambda_{i}}} \exp {-\{{\tr
(\sigma^{-\chi} }X )\}}{ \textbf{1}}_{\mathcal{P}}(X)dX.
\end{eqnarray*}
\end{example}
\begin{example}
In this example, the poset is isomorphic to a rooted tree and
therefore the homogeneous cone ${\cal P}$ is of the $Q_G$ type, that
is, it corresponds to a homogeneous graph (see $[10]$, Theorem 2.2).

\noindent Let the index set be $I=\{1,\ldots,k\}$ and the poset be
defined by$$1\prec i,\;\; i=2,\ldots,k.$$ We can immediately see
here that the undirected graph obtained by dropping the directions
on the directed graphical representation of the poset is the
star-shaped graph with $k$ vertices and undirected edges
$\{(1,i),\;i=2,\ldots,k\}$, which is a homogeneous graph. According
to Theorem 2.2 in $[10]$, the cone ${\cal P}$ is therefore of the
$Q_G$ type. Again, we take $E_{ij}=\reel$ for all appropriate
$(i,j)$. Then $n_{.1}=k-1, \ n_{1.}=0, \ n_1=\frac{k+1}{2}$ and
$n_{.i}=0, \ n_{i.}=1, \ n_i=\frac{3}{2}$, $ \forall i \not = 1$.
The Vinberg algebra ${\mathcal{A}}$ is the vector space of all
$A=(a_{ij}| (i, j) \in I \times I) \in {\mathcal{M}}(I,\reel)$ with
$a_{ij} = 0$ when $i$ and $j$ are not related. We have
$${\mathcal{X}}=\{(\lambda_i,\;i=1,\ldots,k):\; \lambda_1>0, \ \lambda_i>\frac{1}{2},\; \;i=2,\ldots k\},$$
 and the Wishart distribution is given
 by
$$HW_{\chi,\sigma}(dX)=\frac{\pi^{\frac{1}{2}}\prod\limits_{i=1
}^{k}\lambda_{i}^{\lambda_{i}}x_{[1].}^{\lambda_{1}-\frac{k+1}{2}}\prod\limits_{i=2}^{k}
x_{[i].}^{\lambda_{i}-\frac{3}{2}}}{\Gamma(\lambda_{1})\prod
\limits_{i=2
}^{k}\Gamma(\lambda_{i}-\frac{1}{2})\prod\limits_{i=1}^{k}
\sigma_{[i].}^{\lambda_{i}}} \exp {-\{{\tr (\sigma^{-\chi} }X )\}}{
\textbf{1}}_{\mathcal{P}}(X)dX,$$ where  $\chi = \{\lambda_i, \ i\in
I\} \in {\mathcal{X}}$ and $\sigma^{-\chi}$ is defined as in
(\ref{sigma}), for $\sigma \in {\cal{P}}$.

\noindent To illustrate the decomposition of an element of $\cal{P}$
in this example, we consider the case where $k=4$. An element
$X=\pmatrix{x_1&x_{12}&x_{13}&x_{14}\cr x_{21}&x_2&0&0\cr
x_{31}&0&x_3&0\cr x_{41}&0&0&x_4}$ of ${\cal P}$ can then be written
as an incomplete matrix
 $$X^*=\pmatrix{x_1&x_{12}&x_{13}&x_{14}\cr x_{21}&x_2&*&*\cr
x_{31}&*&x_3&*\cr x_{41}&*&*&x_4}.$$

\noindent It can also be decomposed as $X=TDT^\ast$ with $T\in {\cal
T}_l^1, D\in {\cal D}$, where
$$
T^\ast=\pmatrix{1&t_{12}&t_{13}&t_{14}\cr 0&1&0&0\cr 0&0&1&0\cr
0&0&0&1},  \ D=\pmatrix{d_1&0&0&0\cr 0&d_2&0&0\cr 0&0&d_3&0\cr
0&0&0&d_4}.$$ The equation
\begin{eqnarray*}
X &=&\pmatrix{1&0&0&0\cr t_{21}&1&0&0\cr t_{31}&0&1&0\cr
t_{41}&0&0&1} \pmatrix{d_1&0&0&0\cr 0&d_2&0&0\cr 0&0&d_3&0\cr
0&0&0&d_4}\pmatrix{1&t_{12}&t_{13}&t_{14}\cr 0&1&0&0\cr 0&0&1&0\cr
0&0&0&1}\nonumber\\
&=&\pmatrix{d_1&t_{12}d_1&t_{13}d_1&t_{14}d_1\cr
t_{21}d_1&t_{21}d_1t_{12}+d_2&0&0\cr
t_{31}d_1&0&t_{31}d_1t_{13}+d_3&0\cr
t_{41}d_1&0&0&t_{41}d_1t_{14}+d_4}
\end{eqnarray*}
has a unique solution
\begin{eqnarray*}
&&d_1=x_1=x_{[1]\cdot},\;\;t_{12}=\frac{x_{12}}{x_1}\;,t_{13}=\frac{x_{13}}{x_1},\;t_{14}=\frac{x_{14}}{x_1}\\
&&\;d_2=x_2-x_{21}x_1^{-1}x_{12}=x_{[2]\cdot}\;,
d_3=x_3-x_{31}x_1^{-1}x_{13}=x_{[3]\cdot}\;,
d_4=x_4-x_{41}x_1^{-1}x_{14}=x_{[4]\cdot}.
\end{eqnarray*}
Hence, in this example, the cone ${\cal{P}}$ is the set of symmetric
matrices $X$ such that
\begin{eqnarray*}
x_1>0, x_2-x_1^{-1}x_{12}^2>0,\; x_3-x_1^{-1}x_{13}^2>0 \
\textrm{and} \ x_4-x_1^{-1}x_{14}^2>0.
\end{eqnarray*}
\n Let us note that the distribution of $X\in Q_G$ can be obtained
from the distributions of $d_i,i=1,\ldots,k$ and
$t_{1i},i=2,\ldots,k$. This transformation is given explicitly in
$[10]$, Formula (3.15) and Theorem 4.5.\end{example}

Now, we recall that an edge is written as an arrow from its origin
to its destination, and we define the in-degree of a vertex to be
the number of edges having this vertex as their destination. A
vertex is considered a source in a graph if its in-degree is 0 (no
vertices have a source as their destination) and if it is the parent
of at least two vertices. This notion is needed in the following
result.
\begin{prop}
The Vinberg multiplication $TT^{\ast}$ is the same as the standard
matrix multiplication $T.T^{\ast}$ if and only if there is no source
in $I$.
\end{prop}

\begin{Pff} \ $(\Rightarrow)$ Suppose that there exists a source $a$ in $I$. Let $b$ and $c$ be in $I$ such that $a \prec b, \ a \prec c$, and consider
$T=(t_{ij})_{i,j \in I} \in {\cal{T}}_{l}^{+}$ such that $t_{bc}=0$.
Then $T.T^{\ast} \not = TT^{\ast}$,  which contradicts our
assumption and therefore  the source does not exist.

\n $(\Leftarrow)$ Suppose that there exists $T$ such that
$T.T^{\ast} \not = TT^{\ast}$. Then there exist in $I$ two different
and non connected element $i$ and $j$ such that $(TT^{\ast})_{ij}=0$
and $ (T.T^{\ast})_{ij}\not = 0$. As
$(TT^{\ast})_{ij}=\displaystyle\sum_{\begin{array}{c}
  k \preceq i,j \\
\end{array}} t_{ik}t_{jk}$, there exists $k$ in $I$ such that $k \prec i$
and $k \prec j$. This implies that $k$ is a source which is a
contradiction.
\end{Pff}

\n In what follows, we suppose that $I$ has no source and we call
$m$ a maximal element in $I$, if it is not less than any element of
$I$.
\begin{prop}\label{prop}
Let $m$ be a maximal element in $I$ and let $I_{\preceq m}=\{i \in
I, i\preceq m\}$.

\n Then $I_{\preceq m}$ has no source.
\end{prop}

\begin{Pff} \ Suppose that there exists a source $a \in
I_{\preceq m}$ and let $b$ and $c$ in $I_{\preceq m}$ be two
children of $a$ with $a \prec b, \ a \prec c$. As $m$ is a maximal,
we have that $a \not= m$. Also $m \not=b$, because, if $m=b$, then
$c \preceq b$.
 Thus $b\preceq m, \ c \preceq m$ and $a\preceq m$. Hence we obtain $a \prec b, \ a \prec c, \  b \prec m,$

\n $ c \prec m$ and $a \prec m$ which is impossible by Condition (F)
of this section, and the lemma is proved.\end{Pff}

For $T=(t_{ij})_{i,j \in I}$ is in ${{\cal T}}_{l}^{+}$, we define
the element $T_{i\preceq }$ of ${{\cal T}}_{l}$ by $T_{i\preceq }
=(t'_{ij})_{i,j \in I}$, with $t'_{jk}=t_{jk}$ if $i\preceq j, k $
and $t'_{jk}=0$ otherwise. It follows  from  Proposition \ref{prop}
that if $m$ is maximal, then for any $Z \in {\cal{P}}$, we have
$$\det Z_{\preceq m}= \displaystyle \prod_{i\preceq m}Z_{[i].},$$

\n Indeed, $m$  maximal implies that $I_{\preceq m}$ has no source,
and in this case, the Vinberg multiplication is nothing but the
standard multiplication of matrices.

\n Also, we have that for any $i \in I$,
$$\det Z_{\preceq i}= \displaystyle \prod_{i'\preceq i}Z_{[i'].}\ \textrm{and}
\ \det Z_{\prec i}= \displaystyle \prod_{i'\prec i}Z_{[i'].}.$$

\n because $i$ is maximal in $I_{\preceq i}$. Therefore, for all $v
\in I$,
\begin{equation}\label{o2}
Z_{[v].}=\frac{\det Z_{\preceq v}}{\det Z_{\prec v}}\cdot
\end{equation}

We now introduce the notion of a division algorithm in a homogeneous
cone needed for our characterization of the Wishart distribution on
the cone ${\cal{P}}$. This notion is defined in $[6]$ in the case of
a symmetric cone. A division algorithm is a measurable map
\begin{eqnarray}\label{a3}
g:{\mathcal{P}}&\rightarrow&G\nonumber\\
U&\mapsto& g(U)\nonumber
\end{eqnarray}
such that $g(U)(U)=e$. In particular, for $U=TT^{\ast}$, we define
$g(U)=\pi(T^{-1})$, where $\pi$ is defined in (\ref{o6}), then $g$
is a division algorithm, so that if $X=WW^{\ast}\in {{\cal P}}$,
then
$$g(U)(X)=(T^{-1}W)(W^{\ast}(T^{-1})^{\ast}).$$ This algorithm is
the one that we will use in all what follows.
\section{Main characterization result}
In this section, we state and prove our main characterization result
concerning the Wishart distribution on homogeneous cones in the line
of the characterizations given for the ordinary Wishart on symmetric
matrices by Olkin and Rubin $[13]$ and by Casalis and Letac $[6]$.
Our considerations here will be restricted to the case of
homogeneous cones with $E_{ij}=\reel,\; (i,j)\in I\times I.$

It is easy to see that the Laplace transform of the Wishart
distribution $HW_{\chi,\sigma}$ on the homogeneous cone ${\cal{P}}$
is given for $\theta \in {\cal{P}}^{\ast}$ by

$$L_{HW_{\chi},\sigma}(\theta)=\frac{n^{\chi}[(\theta + \sigma^{-\chi})^{\chi}]}{n^{\chi}(\sigma)}.$$
\n

\n As $n^{\chi}(\sigma)=\displaystyle\prod_{i \in
I}\sigma^{\lambda_{i}}_{[i].}n^{\chi}(e),$ this can be written as

$$L_{HW_{\chi},\sigma}(\theta)=\frac{\displaystyle\prod_{i \in
I}[(\theta +\sigma^{-\chi})^{\chi}
]_{[i].}^{\lambda_{i}}}{\displaystyle\prod_{i \in
I}\sigma^{\lambda_{i}}_{[i].}}\cdot$$

\n Recall that for $i \in I$, we denote $I_{\preceq i }=\{j \in I; \
j \preceq i\}$ and $I_{i \preceq }=\{j \in I; \ i \preceq j\}.$ Let
$i_1$ and $i_2$ in $I$, such that $i_1 \not \preceq i_2$ and $i_2
\not \preceq i_1$, we say that $j$ separates $i_{1}$ and $i_{2}$ if
$j \in I_{i_{1}\preceq } \cap I_{i_{2}\preceq }$ and $j\not\in
\{i_{1}, \ i_{2}\}$. In this case, $j$ is called a separator. We
also denote $S_{i}=\{j \in I_{i \preceq }; \ j \ \textrm{is a
separator and} \  \forall \ k \not=j, \ k\not\prec j\}$,
$\mathfrak{S}=\displaystyle\bigcup_{i \in I}S_{i}$ and $S=\{i \in
\mathfrak{S}, \forall \ j \not=i, \ j\not\prec i\}$. This leads to
the following decomposition of an element of ${\cal{P}}$ which will
serve in our characterization result.

\begin{prop}\label{decompos}
Let $Z=TT^{\ast}$ be an element of ${\cal{P}}$ with $T \in {{\cal
T}}_{l}^{+}$. Denote $\wp=\{i \in I, \  I_{ \prec i}=\emptyset\}$
and define, for $i \in I$,
\begin{eqnarray*}\hskip1.5cm Z_{i\preceq }=T_{i\preceq }T_{i\preceq
}^{\ast} \hskip1cm  \textrm{and} \hskip1cm {Z}_{i}= \left \{
\begin{array}{ccccc} Z_{i\preceq }- \displaystyle\sum_{s \in
{S}_{i}} Z_{s\preceq }& \textrm{if} \ & i \in \wp
\hskip1.8cm&&\\Z_{i\preceq }\hskip1.9cm& \textrm{if}\ &i \in S\hskip1.8cm&& \\
{0} \hskip2.1cm&& \textrm{otherwise},\hskip2.3cm && \\
\end{array}\right.
\end{eqnarray*}Then we have that \begin{eqnarray}\label{t1}Z=\displaystyle\sum_{i \in I}Z_{i}=\displaystyle\sum_{i \in
\wp\cup S}Z_{i}.\end{eqnarray}\end{prop}

 \begin{pff} \ In order to prove the equality (\ref{t1}), we compare the blocks of $Z$ and $\displaystyle\sum_{i
\in I}Z_{i}$ on each subalgebra ${\cal{A}}_j=\displaystyle
\prod_{k,l \in I_{j\preceq }}{\cal A}_{kl}$, $j \in I$ which we
denote respectively by $(Z)_j$ and $(\displaystyle\sum_{i \in
\wp\cup S}Z_{i})_j$. From the definition of $(Z_{j\preceq})$, we
first observe that
\begin{eqnarray}\label{Z_j}(Z)_j=(Z_{j\preceq})_j, \
\hskip1cm\forall \ \ j \in I.\end{eqnarray} Now we discuss according
to the position of $j$.

\n If $j \not\in \wp\cup S$, then there exists a unique $i_0 \in
\wp\cup S$ such that $j \in I_{{i_0}\preceq}$. The fact that $j \in
I_{i_0\preceq}$ implies that $(Z_{i_0})_j=(Z_{j\preceq})_j$. From
this and (\ref{Z_j}), we get
$$(Z)_j=(Z_{j\preceq})_j=(Z_{i_0})_j=(\displaystyle\sum_{i \in
\wp\cup S}Z_{i})_j.$$

\n If $j \in S$,  we have $(Z_i)_j=0$, for all $i \in  \wp$.
Therefore
$$(\displaystyle\sum_{i \in \wp\cup S}Z_{i})_j=(\displaystyle\sum_{i \in
S}Z_{i})_j=(Z_{j})_j=(Z_{j\preceq})_j.$$

\n For $j \in \wp$, we have that $(\displaystyle\sum_{i \in
\wp}Z_{i})_j=Z_j=Z_{j\preceq }- \displaystyle\sum_{s \in {S}_{j}}
Z_{s\preceq }$, and we consider separately the cases $S_j=\emptyset$
and $S_j\not=\emptyset$. If $S_j=\emptyset$, then
$(\displaystyle\sum_{i \in S}Z_{i})_j=0$ and
$(Z)_j=(Z_{j\preceq})_j=Z_j=(\displaystyle\sum_{i \in \wp}Z_{i})_j$.

\n If $S_j\not=\emptyset$, $(\displaystyle\sum_{i \in
S}Z_{i})_j=(\displaystyle\sum_{i \in
S_j}Z_{i})_j=(\displaystyle\sum_{i \in S_j}Z_{i\preceq})_j.$
Therefore $(Z)_j=(\displaystyle\sum_{i \in \wp\cup S}Z_{i})_j$.
\end{pff}

\begin{example}
Consider the following poset on $I=\{1,2,3,4\}$ where $1\prec 3, \
1\prec4, \ 2\prec4$. Then $S=S_1=S_2=\{4\}$, $I_{3\preceq}=\{3\}$.
Hence $Z_{3}=(0)$,
$$\small {Z_{4}=Z_{4\preceq }=\left(
         \begin{array}{cccc}
0 &0 &0 &0 \\
0 &0 &0 &0 \\
0 &0 &0 &0 \\
0 &0 &0 & t_{44} \\
         \end{array}  \right)\left(
                      \begin{array}{cccc}
0 &0 &0 &0 \\
0 &0 &0 &0 \\
0 &0 &0 &0 \\
0 &0 &0 &t_{44} \\
                      \end{array}\right)},$$

$$\small{Z_{1}=\left(
\begin{array}{cccc}
t_{11} &0 &0 &0 \\
0 &0 &0 &0 \\
t_{13}&0 &t_{33} &0 \\
t_{14} &0 &0 & t_{44} \\
         \end{array}  \right)\left(
                      \begin{array}{cccc}
t_{11}  &0 &t_{13} &t_{14} \\
0 &0 &0 &0 \\
0 &0 &t_{33}  &0 \\
0 &0 &0 &t_{44} \\
                      \end{array}\right)-Z_{4}},$$
\n and

\n $$\small{Z_{2}=\left(
         \begin{array}{cccc}
0 &0 &0 &0 \\
0 &t_{22} &0 &0 \\
0 &0 &0 &0 \\
0 &t_{24} &0 & t_{44} \\
         \end{array}  \right)\left(
                      \begin{array}{cccc}
0 &0 &0 &0 \\
0 &t_{22} &0 &t_{24}\\
0 &0 &0 &0 \\
0 &0 &0 &t_{44} \\\end{array}\right)- Z_{4}}.$$

\end{example}

\n Consider the Vinberg subalgebra of ${\cal A}$ defined by
${\cal{A}}_{i}=\displaystyle \prod_{k,l \in I_{i\preceq }}{\cal
A}_{kl}$, and denote ${\cal{P}}_{i}$, $G_{i}$ and $e_{i}$
respectively, the corresponding homogeneous cone, connected
component of the identity in $Aut({\cal{P}}_{i}$), and unit element.
Also denote $K= \{g \in G, \ g(e)=e\}$ the orthogonal group of
${\cal{A}}$, and ${K}_{i}=\{k\in K, \ \ k(e_{i})=e_{i} \}$. Finally
let \begin{eqnarray}\label{DA}
g_i:\hskip1.2cm{\mathcal{P}}_i\hskip0.8cm&\rightarrow&G_i\\
U_{i\preceq }=T_{i\preceq }T_{i\preceq }^{\ast}&\mapsto&
g_i(U_{i\preceq })\nonumber
\end{eqnarray}
such that $g_i(U_{i\preceq })(X_{i\preceq })=(T_{i\preceq
}^{-1}W_{i\preceq })(W_{i\preceq }^{\ast}(T_{i\preceq
}^{-1})^{\ast})$, where $X_{i\preceq }=W_{i\preceq }W_{i\preceq
}^{\ast}$. Then it is easy to see that
$$(g(X+Y)X)_{i\preceq }=g_{i}(X_{i\preceq }+Y_{i\preceq })X_{i\preceq }\cdot$$

From now on, a Wishart distribution $HW_{\chi,e^{\chi}}$ will be
called a standard Wishart distribution. Next, we verify that any
Wishart distribution on the homogeneous cone ${\cal{P}}$ may be
standardized by a linear transformation.
\begin{prop}\label{prop1} Let $X$ be a random variable valued in
${\cal{P}}$. Then X is $HW_{\chi,\sigma}$ if and only if there exist
$\rho$ in $G$ such that $\rho(X)$ is $HW_{\chi,e^{\chi}}$.
\end{prop}

\begin{Pff} \ $(\Rightarrow)$ \ Suppose that $X$ is $HW_{\chi,\sigma}$ and write $\sigma=TT^{\ast}$ with $T \in
{\cal{T}}_{l}^{+}$, $\chi=\{\lambda_{i}, \ i\in I\}$. We have that
$\sqrt{\mbox{diag}(\lambda_{i}, \ i \in I)}T^{-1} \in
{\cal{T}}_{l}^{+}$, we then consider the element of $G$,
$\rho=\pi(\sqrt{\mbox{diag}(\lambda_{i}, \ i \in I)}T^{-1})$, where
$\pi$ is defined by (\ref{o6}). Using (\ref{c1}) and (\ref{o7}), we
have for $\theta \in {\cal{P}}^{\ast}$,
\begin{eqnarray}
L_{\rho(X)}(\theta)&=&E(e^{-\mbox{tr}(\theta
\displaystyle\rho(X))})\nonumber\\
&=&E(e^{-\mbox{tr}(\displaystyle\rho^{\ast}(\theta)
X)})\nonumber\\
&=&\frac{\displaystyle\prod _{i \in
I}[(\rho^{\ast}(\theta)+\sigma^{-\chi})^{\chi}]_{[i].}^{\lambda_{i}}}{\displaystyle\prod_{i
\in I}\sigma_{{[i].}}^{\lambda_{i}}}\nonumber\\
&=&\frac{\displaystyle\prod _{i \in
I}[\rho^{-1}(\theta+\rho^{\ast-1}(\sigma^{-\chi}))^{\chi}]_{[i].}^{\lambda_{i}}}{\displaystyle\prod_{i
\in I}\sigma_{{[i].}}^{\lambda_{i}}},\nonumber
\end{eqnarray}
and using the fact that
$n^\chi(g^{-1}(u)(\sigma))=\chi(g^{-1}(u))n^\chi(\sigma)$ and
$\chi(g(u))=\displaystyle\prod_{i \in I} u_{[i].}^{-\lambda_{i}}$
(see $[2]$), we easily verify that
$$\displaystyle \prod_{i \in
I}(g^{-1}(u)\sigma)_{[i].}^{\lambda_i}=\prod_{i \in
I}u_{[i].}^{\lambda_{i}}\sigma_{[i].}^{\lambda_i}\cdot$$ This
applied to $g^{-1}(u)=\rho$ gives
$$L_{\rho(X)}(\theta)=\frac{\displaystyle\prod _{i \in I}[
(\theta+\rho^{\ast-1}(\sigma^{-\chi}))^{\chi}]_{[i].}^{\lambda_{i}}}{\displaystyle\prod_{i
\in I}[\rho(\sigma)]_{{[i].}}^{\lambda_{i}}}\;\nonumber.$$ As
$\rho^{\ast-1}(\sigma^{-\chi})=(\rho(\sigma))^{-\chi}=(e^{\chi})^{-\chi}=e$,
we obtain

$$L_{\rho(X)}(\theta)=\frac{\displaystyle\prod _{i \in I}[
(\theta+e)^{\chi}]_{[i].}^{\lambda_{i}}}{\displaystyle\prod_{i \in
I}(e^{\chi})_{{[i].}}^{\lambda_{i}}},$$ which is the Laplace
transform of a $HW_{\chi,e^{\chi}}$ distribution.

\noindent$(\Leftarrow)$ Suppose that there exist $\rho$ in $G$ such
that $X'=\rho(X)$ is $HW_{\chi,e^{\chi}}$. Then using again
(\ref{c1}) and (\ref{o7}), we have for $\theta \in
{\cal{P}}^{\ast}$,
\begin{eqnarray}
L_{X}(\theta)&=&E(e^{-\mbox{tr}(\theta X)})\nonumber\\
&=&E(e^{-\mbox{tr}(\rho^{\ast-1}(\theta)
X')})\nonumber\\
&=&\frac{\displaystyle\prod _{i \in
I}[(\rho^{\ast-1}(\theta)+e)^{\chi}]_{[i].}^{\lambda_{i}}}{\displaystyle\prod_{i
\in I}(e^{\chi})_{{[i].}}^{\lambda_{i}}}\nonumber\\
&=&\frac{\displaystyle\prod _{i \in I}[
(\theta+\rho^{\ast}(e))^{\chi}]_{[i].}^{\lambda_{i}}}{\displaystyle\prod_{i
\in
I}[(\rho^{\ast}(e))^{\chi}]_{{[i].}}^{\lambda_{i}}}\cdot\nonumber
\end{eqnarray}
\n Hence $X$ is $HW_{\chi,\sigma}$, with $\sigma
=(\rho^{\ast}(e))^{\chi}\in {\mathcal{P}}$.
\end{Pff}

We are now in a position to give our characterization results.
According to Proposition \ref{prop1}, the statements will concern
the Wishart distribution $HW_{\chi,e^{\chi}}$.
\begin{theorem}
\label{theo1} i) Let $X=TT^{\ast}$, with $T=(t_{ij})_{i,j \in I}$ in
${\cal{T}}_{l}^{+},$ be a random variable with Wishart distribution
$HW_{\chi,\sigma}$. Then $\sigma=e^{\chi}$ if and only if the
$t_{ij}, \ i, \ j \in I$ are independent.

\n ii) If $X$ and $Y$ are two independent random variables with
respective Wishart distribution, $HW_{\chi,e^{\chi}}$ and
$HW_{\chi',e^{\chi'}}$, then for $i \in I$, the distribution of
$V_{i\preceq }=g_{i}(X_{i\preceq }+Y_{i\preceq })X_{i\preceq }$ is
${K}_{i}$ invariant, where $g_i$ is defined by (\ref{DA}).
\end{theorem}

Next, we give the reciprocal of this theorem.

\begin{theorem}
\label{theo2} Let $X$ and $Y$ be independent random variables valued
in ${\cal{P}}$. Write $X=TT^{\ast}$ and $Y= MM^{\ast}$, with
$T=(t_{ij})_{i,j \in I}$ and $M=(m_{ij})_{i,j \in I}$ in
${\cal{T}}_{l}^{+}$. Consider the divisions algorithms  $g$ defined
by (\ref{a3}) and $g_i$ defined by (\ref{DA}) and suppose that

\n(i) the $t_{ij}, \ i, \ j \in I$ are independent and the $m_{ij},
\ i, \ j \in I$ are independent,

\n(ii) $X+Y$ is independent of $g(X+Y)(X)$,

\n(iii) for $i \in \wp\cup S$, the distribution of
$V_{i\preceq}=g_{i}(X_{i\preceq }+Y_{i\preceq })X_{i\preceq }$ is
${K}_{i}$ invariant.

\n Then there exist $\chi$, $\chi' \in {\mathcal{X}}$ such that $X
\sim HW_{\chi,e^{\chi}}$ and $Y \sim HW_{\chi',e^{\chi'}}$.

\end{theorem}
Before embarking on the proofs of these theorems, it is worth
mentioning that in the particular case where the Vinberg algebra is
the algebra ${\cal{M}}(I,\reel)$ of $I\times I$ matrices (see
Example \ref{special case}), Theorem \ref{theo1} and Theorem
\ref{theo2} give together the famous Olkin and Rubin
characterization of the ordinary Wishart distribution on symmetric
matrices. In fact in this case we have, $\wp=\{1\}$, $S=\emptyset$,
so that $\wp\cup S=\{1\}$, and it follows that only $K_1=K=\{g \in
G, \ g(e)=e\}$ appears in Theorem the point (iii) of \ref{theo2}. We
also have $e^{\chi}=\lambda e$. Hence Theorem \ref{theo1} becomes:

\n i) Let $X=TT^{\ast}$, with $T=(t_{ij})_{i,j \in I}$ in
${\cal{T}}_{l}^{+},$ be a random variable with Wishart distribution
$W_{\lambda,\sigma}$. Then $\sigma=\lambda e$ if and only if the
$t_{ij}, \ i, \ j \in I$ are independent

\n ii) If $X$ and $Y$ are two independent random variables with
respective Wishart distribution, $W_{\lambda,\lambda e}$ and
$W_{\lambda',\lambda' e}$, then for $i \in I$, the distribution of
$g_{i}(X_{i\preceq }+Y_{i\preceq })X_{i\preceq }$ is ${K}_{i}$
invariant, where $g_i$ is defined by (\ref{DA}), in particular
$g(X+Y)(X)$ is $K$ invariant.
\smallskip

For the proof of the theorems, we need to establish the following
result.
\begin{lemma}\label{lemma1}
Let $X$ and $Y$ be two independent random variables with
respective Wishart distribution, $HW_{\chi,e^{\chi}}$ and
$HW_{\chi',e^{\chi'}}$. For $i \in I$, let $\alpha_{U_{i\preceq }}$ denote the distribution
of $X_{i\preceq }$ conditional on
$X_{i\preceq}+Y_{i\preceq}=U_{i\preceq }$. Then for $f$ in $G_{i}$,
the image measure $f\alpha_{U_{i\preceq }}$ of $\alpha_{U_{i\preceq
}}$ by $f$ is such that $f\alpha_{U_{i\preceq
}}=\alpha_{f(U_{i\preceq })}$.
\end{lemma}

\begin{Pff} \ \  Let $H\;:\; {\cal{A}}\rightarrow {\reel}$ and $F\;:\;
{\cal{A}}\rightarrow {\reel}$ be any continuous functions with
compact support. Let ${\chi}_{i}$ and ${\chi'}_{i}$ be two
multipliers on $G_{i}$ and set
${\chi}^1_{i}={\chi}_{i}{\chi}_{i}'$. Denote $\nu^{{\chi}_{i}},
\ \nu^{{\chi}_{i}'}$ and $\nu^{{\chi}_{i}^{1}}$ the corresponding
equivariant measures concentrated on ${\cal{P}}_{i}$ as defined in
(\ref{meas}) and set $X_{i\preceq }^{1}=f(X_{i\preceq })$ and
$Y_{i\preceq } ^{1}=f(Y_{i\preceq })$. We consider the following
equalities obtained by using the decomposition
 $$\nu^{\chi_{i}}(dX_{i\preceq })\nu^{\chi_{i}'}(dY_{i\preceq })=\alpha_{U_{i\preceq }}(dY_{i\preceq })\nu^{\chi_{i}^1}(dU_{i\preceq })$$ and using the change of variable formula for integrals. For ${U_{i\preceq }}=X_{i\preceq }+Y_{i\preceq }$, we have
\begin{eqnarray} J&=&\int_{{\cal{A}}^{2}}H(U_{i\preceq })F(
X_{i\preceq }^{1})f(\alpha_{U_{i\preceq }})(dX_{i\preceq } ^{1})\nu^{{\chi^1}_{i}}(dU_{i\preceq })\nonumber\\
&=&\int_{{\cal{A}}^{2}}H(U_{i\preceq })(Fof)(X_{i\preceq })\alpha_{U_{i\preceq }}(dX_{i\preceq })\nu^{\chi^1_{i}}(dU_{i\preceq })\nonumber\\
&=&\int_{{\cal{A}}^{2}}H(X_{i\preceq }+Y_{i\preceq
})Fof(X_{i\preceq})\nu^{\chi_{i}}(dX_{i\preceq
})\nu^{\chi'_{i}}(dY_{i\preceq  })\cdot\nonumber
\end{eqnarray}
\n Thus
\begin{eqnarray}
J&=&\int_{{\cal{A}}^{2}}H(f^{-1}(X_{i\preceq }^{1}+Y_{i\preceq } ^{1}))F(X_{i\preceq } ^{1})(f\nu^{\chi_{i}})(dX_{i\preceq }^{1})(f\nu^{\chi'_{i}})(dY_{i\preceq } ^{1})\nonumber\\
&=&\chi_{i}(f^{-1})\int_{{\cal{A}}^{2}}H(f^{-1}(X_{i\preceq }^{1}+Y_{i\preceq } ^{1}))F(X_{i\preceq } ^{1})\nu^{\chi_i}(dX_{i\preceq }^{1})\nu^{\chi'_{i}}(dY_{i\preceq } ^{1})\nonumber\\
&=&\chi_{i}(f^{-1})\int_{{\cal{A}}^{2}}H(f^{-1}(S_{i\preceq
}))F(X_{i\preceq } ^{1})\alpha_{S_{i\preceq }}(dX_{i\preceq }
^{1})\nu^{\chi^1_{i}}(dS_{i\preceq })\;,\nonumber
\end{eqnarray}
where ${S_{i\preceq }}=X^1_{i\preceq }+Y^1_{i\preceq }$. Since
$f^{-1}(S_{i\preceq })=U_{i\preceq }$, we have

\begin{eqnarray}
J&=&\chi^1_{i}(f^{-1})\int_{{\cal{A}}^{2}}H(U_{i\preceq
})F(X_{i\preceq} ^{1})\alpha_{f(U_{i\preceq })}(dX_{i\preceq }
^{1})f^{-1}
\nu^{\chi^1_{i}}(dU_{i\preceq })\nonumber\\
&=&\int_{{\cal{A}}^{2}}H(U_{i\preceq })F(X_{i\preceq }
^{1})\alpha_{f(U_{i\preceq })}(dX_{i\preceq }
^{1})\nu^{\chi^1_{i}}(dU_{i\preceq })\cdot\nonumber
\end{eqnarray}

\noindent Comparing this expression with the definition of $J$, the
lemma is proved.
\end{Pff}

\begin{Pff} {\bf of Theorem \ref{theo1}}\ \ i) $(\Leftarrow)$ See $[2]$
where it is proved that if $X\sim W_{\chi, e^{\chi}}$, then $t_{ij},
\;i,j\in I$ are independent.

\n $(\Rightarrow)$ Suppose that $X=TT^{\ast}$, with $T=(t_{ij})_{i,j
\in I}$ in ${\cal{T}}_{l}^{+}$, is $HW_{\chi,\sigma}$ with
$\sigma\not=e^{\chi}$. We will show that in this case, the $t_{ij}$
are not all independent. As $\sigma\not=e^{\chi}$, there exist
$i\not=j$ such that $\sigma_{ij}\not=0$. Writing $\sigma=WW^{\ast}$,
with $W= (w_{ij})_{i,j \in I}\in {\cal{T}}_{l}^{+}$, we have that
$\sigma_{ij}=\displaystyle\sum_{k \in I}w_{ik}w_{jk}$, and as
$\sigma_{ij}\not=0$, there exists $k \in I$ such that
$w_{ik}w_{jk}\not=0$. From Proposition \ref{prop1}, we have that
$X'=\pi(\sqrt{\mbox{diag}(\lambda_{i}, \ i \in I)} W^{-1}) (X)$ is
$HW_{\chi,e^{\chi}}$. If we set $X'=SS^{\ast}$, with
$S=(s_{ij})_{i,j \in I} \in {\cal{T}}_{l}^{+}$ and
$Z=W\sqrt{\mbox{diag}(\lambda_{i}^{-1}, \ i \in I)}=(z_{ij})_{i,j
\in I}$, we can write
\begin{eqnarray*}
X=\pi^{-1}(\sqrt{\mbox{diag}(\lambda_{i}, \ i \in I)} W^{-1})
(X')=(ZS)(S^{\ast}Z^{\ast})\cdot
\end{eqnarray*}
This implies that $t_{jk}=\displaystyle\sum_{\mu \in
I}z_{j\mu}s_{\mu k}$ for $j \not=k$. Hence
\begin{eqnarray*}
&&E(t_{jk}t_{kk})-E(t_{jk})E(t_{kk})=\displaystyle\sum_{\mu \in
I}z_{j\mu}z_{kk}E(s_{\mu k}s_{kk})-\displaystyle\sum_{\mu \in
I}z_{j\mu}z_{kk}E(s_{\mu k})E(s_{kk})\nonumber\\
&&=\displaystyle\sum_{\mu \not=k}z_{kk}z_{j\mu}E(s_{\mu
k}s_{kk})+z_{jk}z_{kk}E(s_{kk}^{2})-\displaystyle\sum_{\mu \not=k}z_{j\mu}z_{kk}E(s_{\mu k})E(s_{kk})-z_{jk}z_{kk}E(s_{kk})^{2}\cdot\nonumber\\
\end{eqnarray*}
We use the fact that the $s_{ij}$ are independent, because $X'$ is
$HW_{\chi,e^{\chi}}$, to obtain
\begin{eqnarray*}
E(t_{jk}t_{kk})-E(t_{jk})E(t_{kk})&=&z_{jk}z_{kk}(E(s_{kk}^{2})-E(s_{kk})^{2})\nonumber\\
&=&z_{jk}z_{kk}\mbox{var}(s_{kk})\cdot
\end{eqnarray*}
This is different from $0$ because
$z_{kk}=w_{kk}\sqrt{\lambda_{kk}^{-1}}\not=0$,
$z_{jk}=w_{jk}\sqrt{\lambda_{kk}^{-1}}\not=0$ and $s_{kk}$ is not
degenerate, since $s_{kk}^{2}$ is gamma distributed (see $[2]$).

\n ii) Using the notation of Lemma \ref{lemma1} and the fact that
$U$ and $V$ are independent, for two arbitrary continuous functions
with compact support, $H\;:\; {\cal{A}}\rightarrow {\reel}$ and
$F\;:\; {\cal{A}}\rightarrow {\reel}$, we have that
$$E(H(U_{i\preceq })F(V_{i\preceq }))=E(H(U_{i \preceq
})\int_{{\cal{A}}}F(g_i(U_{i\preceq })X_{i\preceq
})\alpha_{U_{i\preceq }}(dX_{i\preceq })).$$

\noindent Writing  $v_{i\preceq }=g_i(U_{i\preceq })X_{i\preceq }$
in the last integral and using Lemma 4.1, we get
$$E(H(U_{i\preceq })F(V_{i\preceq }))=E(H(U_{i\preceq
})\displaystyle\int_{{\cal{A}}}F(V_{i\preceq
})\alpha_{e_{i}}(dV_{i\preceq }))$$

\n This proves (ii).
\end{Pff}

\begin{Pff} {\bf of Theorem \ref{theo2}}\ \ Let $i \in I$ such that $X_i \not=0$, where $X_i$ is defined as in (\ref{t1}) and let $
U=X+Y=WW^{\ast}$, $V=g(X+Y)(X)$. Without loss of generality, we can
suppose that $S_{i}$ has just one element; $S_i=\{s\}$. We will
consider first the case where $i\not\in S$. In this case, from the
hypotheses of independence, we have, for $A_{1}, \ A_{2}, \ B_{1}, \
B_{2}$ and $C_{1}$  in ${\cal{P}}^{\ast}$,
\begin{eqnarray}\label{t4}
&&E(\exp{{\mathrm{tr}}(A_{1}W_{i\preceq }-A_{2}W_{s\preceq
}+B_{1}W_{i\preceq }W_{i\preceq }^{\ast}-B_{2}W_{s\preceq
}W_{s\preceq
}^{\ast}+C_{1}V_{i\preceq })})\nonumber\\
&&=E(\exp{{\mathrm{tr}}(A_{1}W_{i\preceq }-A_{2}W_{s\preceq
}+B_{1}W_{i\preceq }W_{i\preceq }^{\ast}-B_{2}W_{s\preceq
}W_{s\preceq }^{\ast})})E(\exp{{\mathrm{tr}}( C_{1}V_{i\preceq })})\nonumber\\
&&=f_{i}(A_{1},A_{2},B_{1},B_{2})h_{i}(C_{1}),
\end{eqnarray}
\n where
$f_{i}(A_{1},A_{2},B_{1},B_{2})=E(\exp{{\mathrm{tr}}(A_{1}W_{i\preceq
}-A_{2}W_{s\preceq }+B_{1}W_{i\preceq }W_{i\preceq
}^{\ast}-B_{2}W_{s\preceq }W_{s\preceq }^{\ast})})$

\n and $h_{i}(C_{1})=E(\exp{{\mathrm{tr}} C_{1}V_{i\preceq }})$.

\n We adopt the notations:

$$(f_{i})^{jk}=\frac{\partial f_{i}}{\partial
(A_{1})_{jk}}, \ (f_{i})_{jk}=\frac{\partial f_{i}}{\partial
(B_{1})_{jk}} \ \mbox{and} \ (h_{i})_{jk}=\frac{\partial
h_{i}}{\partial (C_{1})_{jk}},$$

\n and we define, for $B \in {\cal{P}}^{\ast}$,
${\widetilde{f_{i}}}(B)=E(\exp{\textrm{tr} B(X+Y)_{i}})\cdot$

\n It is clear that if $A_{1}=A_{2}=0$ and $B_{1}=B_{2}=B$, then
\begin{equation}\label{o8}
(\widetilde{f}_{i})_{jk}=(f_{i})_{jk} \hskip1cm \textrm{when}
\hskip1cm \{j,k\} \not \subset S_i
\end{equation}
\n and
\begin{equation}\label{o9}
(\widetilde{f}_{i})_{jk,ln}=(f_{i})_{jk,ln} \hskip1cm \textrm{when}
\hskip1cm \{j,k\} \not \subset S_i \ \textrm{and} \ \{l,n\} \not
\subset S_i.
\end{equation}

Since we can differentiate under the expectation, there is a
relation between the second partial derivatives with respect to
$A_{1}$ and the first partial derivatives with respect to $B_{1}$,
namely we have
\begin{equation}\label{o1}
\displaystyle\sum_{\alpha}\frac{\partial^{2}}{\partial
(A_{1})_{j\alpha}\partial(A_{1})_{k\alpha}}\frac{\partial^{q}f_{i}}{\partial
t_{1}\cdot\cdot\cdot\partial
t_{q}}=\frac{\partial}{\partial(B_{1})_{jk}}\frac{\partial^{q}f_{i}}{\partial
t_{1}\cdot\cdot\cdot\partial t_{q}},
\end{equation}

\n where $t_{1},\cdot\cdot\cdot,t_{q}$ are any arguments of $f_{i}$.

Differentiating (\ref{t4}) successively with respect to $(A_{1})_{j
\lambda}, \ (C_{1})_{\lambda \mu}, \ \mbox{and} \ (A_{1})_{k \mu}$,
for $\{j,k\} \not\subset~S_i $ and summing over $\lambda$ and $\mu$,
give the basic differential equations
\begin{eqnarray}\label{t5}
&&E((X_{i})_{jk}\exp{{\mathrm{tr}}(A_{1}W_{i\preceq
}-A_{2}W_{s\preceq }+B_{1}W_{i\preceq }W_{i\preceq
}^{\ast}-B_{2}W_{s\preceq
}W_{s\preceq }^{\ast}+C_{1}V_{i\preceq })})\nonumber\\
&&=\displaystyle\sum_{\lambda,\mu \in I_{i\preceq
}}f_{i}^{j\lambda,k\mu}(h_{i})_{\lambda\mu}.
\end{eqnarray}
\n And differentiating (\ref{t5}) successively with respect to
$(A_{1})_{l \nu}, \ (C_{1})_{\nu \sigma} \ \mbox{and} \ (A_{1})_{n
\sigma}, \ \textrm{for} \ \{l,n\} \not \subset S_i $ and summing, we
obtain
\begin{eqnarray}\label{t6}
&&E((X_{i})_{jk}(X_{i})_{ln}\exp{{\mathrm{tr}}(A_{1}W_{i\preceq
}-A_{2}W_{s\preceq }+B_{1}W_{i\preceq }W_{i\preceq
}^{\ast}-B_{2}W_{s\preceq }W_{s\preceq }^{\ast}+C_{1}V_{i\preceq
})})\nonumber\\
&&=\displaystyle\sum_{\lambda,\mu,\nu,\sigma \in I_{i\preceq
}}f_{i}^{j\lambda,k\mu,l\nu,n\sigma}(h_{i})_{\lambda\mu,\nu\sigma}\cdot
\end{eqnarray}
Now, we use the hypothesis of invariance of the distribution of
$V_{i\preceq}$ by the orthogonal group $K_i$. In terms of Laplace
transform, we have that for any $k \in K_i$,
\begin{eqnarray}\label{inv}h_{i}(C_{1})=E(\exp{{\mathrm{tr}}C_{1}V_{i\preceq
}})=E(\exp{{\mathrm{tr}}C_{1}k(V_{i\preceq })}).
\end{eqnarray}
By the choice of suitable $k$ in $K_i$ and differentiation of
(\ref{inv}), we first establish that there exists a real constant
$\theta_i$ such that
\begin{eqnarray}\label{t2}
(h_{i})_{jm}(0)&=&\frac{\partial h_{i}(C_{1})}{\partial
C_{jm}}\big|_{C_{1}=0}=\theta_i\delta_{jm}, \  \ \mbox{for}\ j ,\ m
\in I_{i\preceq},
\end{eqnarray}
where $\delta_{jm}$ is the Kronecker delta, then that

\n$(h_i)_{jj,jm}(0,0)=(h_i)_{jj,ml}(0,0)=(h_i)_{jm,jl}(0,0)=(h_i)_{jm,ln}(0,0)=0$,
for all different $j,\ m, \ l$  and $n$ in $I_{i\preceq}$. Now, we
set $k=(t_{jm})_{ j,m \in I_{i\preceq}}$ in (\ref{inv}) and we
differentiate to get
$$(h_{i})_{jm,ln}(0,0)=\displaystyle\sum_{\alpha,\beta,\eta,\delta}t_{j\alpha}t_{m\beta}t_{l\eta}t_{n\delta}(h_{i})_{\alpha\beta,\eta\delta},$$
which yields
\begin{eqnarray}\label{t3}
(h_{i})_{jm,ln}(0,0)&=&\eta_{i}\delta_{jm}\delta_{ln}+\xi_{i}[\delta_{jl}\delta_{mn}+\delta_{jn}\delta_{ml}],
\ \ j,\ m, \ l,\ n \in I_{i\preceq},
\end{eqnarray}
\n where $ \eta_{i} \ \mbox{and} \ \xi_{i}$ are real constants.

We set in (\ref{t5}) and (\ref{t6}) $A_{1}=A_{2}=C_{1}=0$ and
$B_{1}=B_{2}=B$. Taking into account (\ref{o1}), (\ref{t2}) and
(\ref{t3}), we obtain, for $\{j,k\} \not \subset S_i \ \textrm{and}
\ \{l,n\} \not \subset S_i$,
\begin{eqnarray}\label{t7} E((X_{i})_{jk}\exp{{\mathrm{tr}}
B(X+Y)_{i}})&=&\displaystyle\sum_{\lambda, \mu \in
I_{i\preceq }}(\widetilde{f}_{i})^{j\lambda,k\mu}(h_{i})_{\lambda\mu} \nonumber\\
&=&\theta_i\displaystyle\sum_{\lambda, \mu \in
I_{i\preceq }}(\widetilde{f}_{i})^{j\lambda,k\mu}\delta_{\lambda\mu} \nonumber\\
&=&\theta_i(\widetilde{f}_{i})_{jk},
\end{eqnarray}
and
\begin{eqnarray}\label{t8}
E((X_{i})_{jk}(X_{i})_{ln}\exp{{\mathrm{tr}}
B(X+Y)_{i}})&=&\displaystyle\sum_{\lambda,\mu,\nu,\sigma \in
I_{i\preceq }}f_{i}^{j\lambda,k\mu,l\nu,n\sigma}[\eta_{i}\delta_{\lambda\mu}\delta_{\nu\sigma}+\xi_{i}(\delta_{\lambda\nu}\delta_{\mu\sigma}+\delta_{\lambda\sigma}\delta_{\mu\nu})]\nonumber\\
&=&\eta_{i}(\widetilde{f}_{i})_{jk,ln}+\xi_{i}[(\widetilde{f}_{i})_{jl,kn}+(\widetilde{f}_{i})_{jn,lk}]\cdot
\end{eqnarray}
\n In the case where $i$ is a separator, the system (\ref{t4}) is
replaced by
\begin{eqnarray*} E(\exp{{\mathrm{tr}}(A_{1}W_{i\preceq
}+B_{1}W_{i\preceq }W_{i\preceq }^{\ast}+C_{1}V_{i\preceq
}}))&=&E(\exp{{\mathrm{tr}}(A_{1}W_{i\preceq }+B_{1}W_{i\preceq
}W_{i\preceq }^{\ast})})E(\exp{{\mathrm{tr}} C_{1}V_{i\preceq }})\nonumber\\
&=&f_{i}(A_{1},B_{1})h_{i}(C_{1}).
\end{eqnarray*}
and by the same reasoning, we also obtain equations (\ref{t7}) and
(\ref{t8}), but for $j, \ k, \ l, \ n$ in $I_{i\preceq}$.

\n Now, to solve (\ref{t7}) and (\ref{t8}), we introduce the
functions
\begin{equation}\label{o3}
\varphi_{i}(B)=E(\exp {\textrm{tr} BX_{i}}), \ \ \psi_{i}(B)=E(\exp
{\textrm{tr} BY_{i}}).
\end{equation}
\n Then (\ref{t7}) and (\ref{t8}) can be written as
\begin{equation}\label{t9}
(\varphi_{i})_{jk}\psi_{i}=\theta_i(\varphi_{i}\psi_{i})_{jk},
\end{equation}
\begin{equation}\label{t10}
(\varphi_{i})_{jk,ln}\psi_{i}=\eta_{i}(\varphi_{i}\psi_{i})_{jk,ln}+\xi_{i}[(\varphi_{i}\psi_{i})_{jl,kn}+(\varphi_{i}\psi_{i})_{jn,lk}]\cdot
\end{equation}
 Equation (\ref{t9}) implies in particular that
\begin{equation}\label{t11}
\varphi_{i}=(\varphi_{i}\psi_{i})^{\theta_i}.
\end{equation}

\n Let $\phi_{i}$ be such that, for $B$ sufficiently close to zero,
\begin{eqnarray}\label{phi}
\varphi_{i}\psi_{i}=\exp \phi_{i}.
\end{eqnarray}
\n Then (\ref{t11}) and (\ref{phi}) lead to the differential
equation system defined for $j, \ k, \, l, \ n \in I_{i\preceq}$ by
\begin{eqnarray}\label{t13}
\theta_i(\phi_{i})_{jk,ln}+\theta_i^{2}(\phi_{i})_{jk}(\phi_{i})_{ln}
&=&\xi_{i}[(\phi_{i})_{jl,kn}+(\phi_{i})_{jl}(\phi_{i})_{kn}
+(\phi_{i})_{jn,lk}+(\phi_{i})_{jn}(\phi_{i})_{lk}]\nonumber\\
&+&\eta_i[(\phi_{i})_{jk,ln}+(\phi_{i})_{jk}(\phi_{i})_{ln}].
\end{eqnarray}
\n The general solution of this system is
\begin{eqnarray}\label{t15}
\phi_{i}(B)&=&\beta_i\ln(B+D)^{\chi_{1}}_{[i].}+a_i,\nonumber
\end{eqnarray}
\n where $D \in {\cal{P}}^{\ast}, \chi_{1} \in {\mathcal{X}}$, and
$a_i \in \reel$ are arbitrary constants and $\beta_i$ is a constant
depending on $\theta_i, \ \xi_{i}, \ \eta_{i}$. This is in
particular, justified by the fact that
\begin{eqnarray*}
(B+D)^{\chi_{1}}_{[i].}&=&(B+D_{i})^{\chi_{1}}_{[i].},\nonumber
\end{eqnarray*}
\n where $D_i \in {\cal{P}}_{i}^{\ast}.$  Now, if $A=R^{\ast}R \in
{\cal{P}}^{\ast}$, we denote
$A_{\preceq^{\opp}}=R^{-1}(R^{-1})^{\ast} \in {\cal{P}}$. Therefore,
we use (\ref{o2}) and a result in $[2]$ due to Andersson and Wojnar
which says that, for $\theta=Z^{\ast}Z$, where $Z\in {\cal
T}_l^{+}$, $\theta^{\chi}=Z^{-1}e^{\chi}(Z^{\ast})^{-1}$,
 to write
\begin{eqnarray*}
\phi_{i}(B)=\beta_i\ln((B+D)_{\displaystyle{\preceq}^{\opp}})_{[i].}+a'_i=\beta_i[\ln
\det (B+D)_{\displaystyle{\preceq}^{\opp}i}-\ln\det
(B+D)_{\displaystyle{\prec}^{\opp}i}]+ a'_i,
\end{eqnarray*}

\n where $a'_i$ is a constant.

\n Thus
$$\varphi_{i}(B)=\frac{((B+D)_{\displaystyle{\preceq}^{\opp}i})_{[i].}^{\beta_{i}\theta_{i}}}{(D_{\displaystyle{\preceq}^{\opp}i})_{[i].}^{\beta_{i}\theta_{i}}}
\ \ \textrm{and} \ \
\psi_{i}(B)=\frac{(B+D)\displaystyle{\preceq}^{\opp}i)_{[i].}^{\beta_{i}(1-\theta_{i})}}{(D_{\displaystyle{\preceq}^{\opp}i})_{[i].}^{\beta_{i}(1-\theta_{i})}}\cdot$$

\n We now observe that the fact that the $t_{ij}, \ i, \ j \in I$
are independent implies that the $X_{i}, \ i \in I$ are independent.
This is important for the calculation of the Laplace transform of
$X=\displaystyle\sum_{i \in I}X_i$. In fact, for each $i \in I$, by
its very definition, $X_i$ depends only on the $t_{ij}$ such that
$i\preceq j$ and on the $t_{kj}$ such that $j,k \not\in S_i \
\textrm{and} \ i \preceq k,j$, which are different from the ones on
which any other component $X_{i'}$, $i'\not=i$ depends. In other
words, there exists a partition $(\tau_i)_{i\in I}$ of the set
$\{t_{jk}, \ j, \ k \in I\}$ such that each $X_i$ depends only on
the $t_{kj}$ in $\tau_i$. Similarly the independence of the $a_{ij},
\ i, \ j \in I$ implies the independence of the $Y_i, \ i \in I$.
Finally denoting $\chi=\{\beta_{i}\theta_{i}, \ i\in I \}$ and
$\chi'=\{\beta_{i}(1-\theta_{i}), \ i\in I \}$,  we obtain the
Laplace transforms $\varphi$ of $X=\displaystyle\sum_{i \in I}X_i$
and $\psi$ of $Y=\displaystyle\sum_{i \in I}Y_i$ as
$$\varphi(B)=\displaystyle\prod _{i \in I}\varphi_{i}(B)=\frac{\displaystyle\prod _{i \in I}[(B+D)^{\chi}]_{[i].}^{\beta_{i}\theta_{i}}}{\displaystyle\prod_{i \in I}[D^{\chi}]_{{[i].}}^{\beta_{i}\theta_{i}}}, \hskip1.5cm \psi(B)=\displaystyle\prod _{i \in I}\psi_{i}(B)=\frac{\displaystyle\prod _{i \in I}[(B+D)^{\chi'}]_{[i].}^{\beta_{i}(1-\theta_{i})}}{{\displaystyle\prod_{i \in
I}[D^{\chi'}]_{[i].}^{\beta_{i}(1-\theta_{i})}}}.$$

\n Thus $X\sim HW_{\chi,\sigma} \ \textrm{and} \  Y \sim
HW_{\chi',\sigma'}$, where $\sigma=D^{\chi}$ and
$\sigma'=D^{\chi'}$. Invoking Theorem \ref{theo1} i), we have
necessarily $D=e$. This concludes the proof of the Theorem
\ref{theo2}.
\end{Pff}

\end{document}